\documentclass[12pt]{amsart} 
\usepackage[english]{babel}
\usepackage[utf8]{inputenc}
\usepackage{amssymb}
\usepackage{amsmath}
\usepackage{amscd}
\usepackage{latexsym}
\usepackage{amsthm}
\usepackage{mathrsfs}
\usepackage{mathabx}
\usepackage{graphicx}
\usepackage{multicol}
\usepackage{footmisc}
\usepackage{color}
\theoremstyle{plain}
\newtheorem{theorem}{Theorem}[section]

\newtheorem{proposition}[theorem]{Proposition}

\newtheorem{definition}[theorem]{Definition}
\newtheorem{remark}[theorem]{Remark}
\newtheorem{example}[theorem]{Example}

\begin{document}

\title{Examples of irreducible symplectic varieties}
\author{Arvid Perego}

\begin{abstract}
Irreducible symplectic manifolds are one of the three building blocks of compact K\"ahler manifolds with numerically trivial canonical bundle by the Beauville-Bogomolov decomposition theorem. There are several singular analogues of irreducible symplectic manifolds, in particular in the context of compact K\"ahler orbifolds, and in the context of normal projective varieties with canonical singularities. In this paper we will collect their definitions, analyze their mutual relations and provide a list of known examples.
\end{abstract}

\maketitle

\section{Introduction}

A central problem in complex geometry is the classification of Ricci-flat compact K\"ahler manifolds. By Yau's theorem \cite{Yau}, these are exactly the compact K\"ahler manifolds whose first Chern class is zero in $H^{2}(X,\mathbb{R})$ or, equivalently, whose canonical bundle is numerically trivial. This implies that the Kodaira dimension is zero.

In dimension 1, compact K\"ahler manifolds of Kodaira dimension zero are exactly elliptic curves. The birational classification of compact complex surfaces shows that compact K\"ahler surfaces of Kodaira dimension zero are K3 surfaces, $2-$dimensional complex tori, Enriques surfaces and bielliptic surfaces.

Among compact K\"ahler manifolds of Kodaira dimension zero, a very special role is played by those manifolds whose canonical bundle is trivial (which are sometimes called \textit{Calabi-Yau manifolds}).

The first family of examples is given by \textit{complex tori}, i. e. quotients of a complex vector space $V$ of dimension $n$ by a rank $2n$ lattice $\Gamma$ in $V$. A projective complex torus is called \textit{Abelian variety}: in dimension 1 all complex tori are projective, and they are exactly the elliptic curves. In higher dimension there are complex tori which are not projective. In any case, if $X$ is a complex torus of dimension $n$, then $\pi_{1}(X)\simeq\mathbb{Z}^{2n}$.

A second family of examples is given by \textit{special unitary manifolds}, i. e. compact K\"ahler manifolds with trivial canonical bundle, whose every finite \'etale covering has no non-trivial holomorphic $p-$form for $0<p<n$ (where $n$ is the complex dimension). They are projective for $n\geq 3$, and have finite fundamental group, which is trivial if $n$ is even (see Proposition 2 of \cite{B}). Simply connected special unitary manifolds are called \textit{irreducible Calabi-Yau manifolds}.

A third family of examples is given by irreducible symplectic manifolds, that will be our main interest. We recall that if $X$ is a complex manifold, a \textit{holomorphic symplectic form} on $X$ is a closed, holomorphic $2-$form $\sigma$ on $X$ which is everywhere non-degenerate. A \textit{holomorphic symplectic manifold} is a complex manifold admitting a holomorphic symplectic form. Among holomorphic symplectic manifolds we find all even-dimensional complex tori.

If $X$ is a holomorphic symplectic manifold, then its complex dimension is even. If $2n$ is the complex dimension of $X$, and $\sigma$ is a holomorphic symplectic form on $X$, then $\sigma^{n}$ is a nowhere vanishing section of $K_{X}$: it follows that a holomorphic symplectic manifold has trivial canonical bundle. Moreover, as $\sigma$ is closed, then it defines a nontrivial cohomology class in $H^{0}(X,\Omega_{X}^{2})$. 

\begin{definition}
An \textit{irreducible symplectic manifold} is a compact, K\"ahler, holomorphic symplectic manifold $X$ which is simply connected and such that $h^{2,0}(X)=1$.
\end{definition}

\begin{remark}
As shown in Propositions 3 and 4 of \cite{B}, an irreducible symplectic manifold $X$ of dimension $n$, then for every $0\leq p\leq n$ we have 
\begin{equation}
\label{eq:hodgenum}
h^{p,0}(X)=\left\{\begin{array}{ll} 0, & {\rm if}\,\,p\,\,{\rm is}\,\,{\rm odd}\\ 1, & {\rm otherwise}\end{array}\right.
\end{equation}
Conversely, by Proposition A.1 of \cite{HNW} every compact, K\"ahler, holomorphic symplectic manifold $X$ of complex dimension $n$ such that $h^{p,0}(X)$ is as in (\ref{eq:hodgenum}) is simply connected, and hence an irreducible symplectic manifold.
\end{remark}

\begin{remark}
\label{remark:holsp}
{\rm As shown in Proposition 4 of \cite{B}, a compact K\"ahler manifold of complex dimension $2n$ is a holomorphic symplectic manifold if and only if its holonomy group is contained in the symplectic group $Sp(r)$. The holonomy group is precisely $Sp(r)$ if and only if $X$ is an irreducible symplectic manifold.}
\end{remark}

As already recalled, Ricci-flat compact K\"ahler surfaces are K3 surfaces, complex tori of dimension 2, Enriques surfaces and bielliptic surfaces. Among them, K3 surfaces are both irreducible Calabi-Yau and irreducible symplectic. 

An Enriques surface is a finite quotient of a K3 surface by a fixed point free involution, while a bielliptic surface is quotient of an Abelian surface (product of two elliptic curves) by the free action of a finite Abelian group. It follows that a Ricci-flat compact K\"ahler surface has a finite \'etale covering wich is either an irreducible symplectic (Calabi-Yau) surface or a complex torus of dimension 2.

The same phenomenon occurs for higher dimensional Ricci-flat compact K\"ahler manifolds, giving to complex tori, irreducible Calabi-Yau manifolds and irreducible symplectic manifolds a special role in the classification. This is the content of the following, which goes under the name of \textit{Beauville-Bogomolov Decomposition Theorem}, and whose proof (based on several results of differential geometry) is contained in \cite{Bog2} and \cite{B}.

\begin{theorem}
\label{theorem:bbdecomp}
Let $X$ be a Ricci-flat compact K\"ahler manifold. Then $X$ has a finite \'etale covering $f:Y\longrightarrow X$, where $Y$ is a product of complex tori, irreducible Calabi-Yau manifolds and irreducible symplectic manifolds.
\end{theorem}

While it is not difficult to provide examples of complex tori and irreducible Calabi-Yau manifolds, it is a hard problem to give examples of irreducible symplectic manifolds. The list of all the known deformation classes is very short:
\begin{enumerate}
 \item Irreducible symplectic surfaces are exactly K3 surfaces.
 \item For $n\geq 2$, the Hilbert scheme $Hilb^{n}(S)$ of $n$ points on a K3 surface $S$ is an irreducible symplectic manifold of dimension $2n$ (see Th\'eor\`eme 3 of \cite{B}). 
 \item For $n\geq 2$, let $T$ be a complex torus of dimension 2 and $s:Hilb^{n+1}(T)\longrightarrow T$ the sum morphism. Then $Kum^{n}(T):=s^{-1}(0_{S})$, called \textit{generalized Kummer variety}, is an irreducible symplectic manifold of dimension $2n$ (see Th\'eor\`eme 4 of \cite{B}).
 \item Two more deformation examples, $OG_{6}$ of dimension 6 and $OG_{10}$ of dimension 10, were constructed by O'Grady (in \cite{OG3} and \cite{OG2} respectively) as a symplectic resolution of some singular moduli spaces of semistable sheaves on a projective K3 or on an Abelian surface. 
\end{enumerate}

The previous examples form different deformation classes because even if they can have the same dimension, they have different second Betti number (which is, following the previous ordering, 22, 23, 7, 8 and 24, see \cite{B}, \cite{OG2}, \cite{OG3} and \cite{R}).

\begin{remark}
{\rm Different constructions of examples of irreducible symplectic manifolds were presented in \cite{BD}, \cite{IR}, \cite{OG4}, \cite{DV} (deformation equivalent to $Hilb^{2}(K3)$) and \cite{LLSvS} (deformation equivalent to $Hilb^{4}(K3)$). Moduli spaces of stable sheaves or of Bridgeland stable complexes on projective K3 surfaces or on Abelian surfaces give rise to examples of irreducible symplectic manifolds which are deformation equivalent to Hilbert schemes of points on K3 surfaces or to generalized Kummer varieties on Abelian surfaces.}   
\end{remark}

Theorem \ref{theorem:bbdecomp} is stated only for compact K\"ahler manifolds. Anyway, if $X$ is a complex projective manifold with Kodaira dimension zero, the minimal model $Y$ of $X$ (whose conjectural existence is predicted by the Minimal Model Program) is a projective variety which is birational to $X$, and has terminal singularities and nef canonical divisor. Assuming the abundance conjecture, it follows that a multiple of the canonical divisor $K_{Y}$ of $Y$ is trivial. 

For the classification of projective varieties whose Kodaira dimension is 0 it is then central to extend Theorem \ref{theorem:bbdecomp} to normal projective varieties having terminal singularities and torsion (i. e. numerically trivial by Theorem 8.2 of \cite{Kaw}) canonical divisor. This implies the need for a definition of singular analogues of irreducible Calabi-Yau and symplectic manifold.

Various definitions have been proposed and studied over the years, and the main purpose of this survey is to present a list of definitions of irreducible Calabi-Yau and symplectic varieties which can be found in the literature, together with a list of known example.

\section{Irreducible symplectic varieties}

There are two main settings we will consider: orbifolds and varieties with canonical singularities. The second one is more natural for the purposes of the Minimal Model Program, while the first one has the advantage to be more similar to the smooth case (and it was the first one to be considered).

\subsection{Orbifolds}

A first generalization of the decomposition theorem was obtained by Campana in \cite{Ca}, in the setting of orbifolds. We refer the reader to \cite{DK} or to \cite{Fu} (where orbifolds are called \textit{V-manifolds}) for precise definitions and results on these analytic spaces.

An \textit{orbifold} is a connected, para-compact analytic space $X$ such that for every $x\in X$ there is an open neighborhood $U$ of $x$ in $X$, an open subset $V$ of $\mathbb{C}^{n}$ (with respect to the Euclidean topology) and a finite group $G$ of automorphisms of $V$ such that there is an isomorphism $\phi:V/G\longrightarrow U$. The composition $\pi:V\longrightarrow U$ of the projection $V\longrightarrow V/G$ with $\phi$ is called \textit{uniformization map}, and the triple $(V,G,\pi)$ a \textit{uniformizing system} for $U$.

It follows from Proposition 1.3 of \cite{Bla} and by Proposition 5.15 of \cite{KM2} that an orbifold is normal, $\mathbb{Q}-$factorial, Cohen-Macaulay and has rational singularities.

On an orbifold $X$ of dimension $n$ one can then always find a uniformizing open cover $\{U_{i}\}_{i\in I}$, i. e. for each $i\in I$ there is a uniformizing system $(V_{i},G_{i},\pi_{i})$ for $U_{i}$. The sheaf of differential $p-$forms on $X$ is denoted $\mathcal{A}^{p}_{X}$, and by definition it is the sheaf which restricted to $U_{i}$ is $\pi_{i*}(\mathcal{A}^{p}_{V_{i}})^{G_{i}}$, i. e. the push-forward under $\pi_{i}$ of the $G_{i}-$invariant part of the sheaf of differential $p-$forms on $V_{i}$. Similarly one defines the sheaf of holomorphic $p-$forms on $X$, denoted $\Omega_{X}^{p}$.

A differential $2-$form $\omega\in\mathcal{A}^{2}_{X}(X)$ is a \textit{K\"ahler form} for $X$ if for every $x\in X$ and every uniformizing system $(V,G,\pi)$ of an open neighborhood $U$ of $x$, we have that $\omega_{|U}\in(\mathcal{A}^{2}_{V})^{G}(V)$ is a K\"ahler metric on $V$.

Moreover, one can still define the notion of first Chern class for orbifolds: the canonical sheaf $K_{X}$ (which, restricted to every uniformized open subset $U$ with uniformizing system $(V,G,\pi)$, is just $\pi_{*}K_{V}$) is not in general locally free, but $K_{X}^{m}$ is if $m$ is a multiple of $|G|$. It follows that if $X$ is compact, then there is $m\gg 0$ such that $K_{X}^{m}$ is locally free, and one defines $c_{1}(X):=-\frac{1}{m}c_{1}(K_{X}^{m})$.

As shown by Theorem 1.1 of \cite{Ca}, if $X$ is a compact, connected K\"ahler orbifold such that $c_{1}(X)=0$ in $H^{2}(X,\mathbb{R})$, every K\"ahler class is represented by a unique Ricci-flat K\"ahler metric on $X$.

\subsubsection{The decomposition theorem for orbifolds}

As in the smooth case, one has two important kinds of Ricci-flat compact K\"ahler orbifolds which arise. The first one generalizes irreducible Calabi-Yau manifolds.

\begin{definition}
An \textit{irreducible Calabi-Yau orbifold} is a compact, connected K\"ahler orbifold $X$ having simply connected smooth locus, $K_{X}\simeq\mathcal{O}_{X}$ and $H^{0}(X,\Omega_{X}^{p})=0$ for every $0<p<\dim(X)$.
\end{definition}

By Proposition 6.6 of \cite{Ca} an irreducible Calabi-Yau orbifold of dimension $n$ is equivalently a compact, connected K\"ahler orbifold whose smooth locus is simply connected, and which has a Ricci-flat K\"ahler metric whose holonomy group is $SU(n)$.

The second kind of Ricci-flat compact K\"ahler orbifolds is a generalization of irreducible symplectic manifolds. A closed holomorphic $2-$form on $X$, i. e. a global section of $\Omega_{X}^{2}$, which is everywhere non-degenerate is called \textit{symplectic form}. An orbifold admitting a holomorphic symplectic form is called \textit{holomorphic symplectic orbifold}. As for holomorphic symplectic manifolds, a holomorphic symplectic orbifold has even complex dimension and trivial canonical sheaf.

\begin{definition}
An \textit{irreducible symplectic orbifold} is a compact, connected K\"ahler orbifold which is holomorphic symplectic, has simply connected smooth locus and $h^{0}(X,\Omega_{X}^{p})=1$.
\end{definition}

By Proposition 6.6 of \cite{Ca} an irreducible symplectic orbifold of dimension $2n$ is equivalently a compact, connected K\"ahler orbifold whose smooth locus is simply connected, and which has a Ricci-flat K\"ahler metric whose holonomy group is $Sp(n)$.

The decomposition theorem for Ricci-flat compact K\"ahler orbifolds is the following (see Theorem 6.4 of \cite{Ca}):

\begin{theorem}
\label{theorem:bbdecomporb}
Let $X$ be a Ricci-flat compact K\"ahler orbifold. Then $X$ has a finite quasi-\'etale covering $f:Y\longrightarrow X$, where $Y$ is a product of complex tori, irreducible Calabi-Yau orbifolds and irreducible symplectic orbifolds.
\end{theorem}

A finite quasi-\'etale morphism is a finite morphism which is \'etale in codimension 1. The statement of Theorem 6.4 of \cite{Ca} is more precise about the finite quasi-\'etale covering: it is indeed a finite orbifold covering (see D\'efinition 5.1 of \cite{Ca}).

The proof of Theorem \ref{theorem:bbdecomporb} relies first on a generalization to orbifolds of the de Rham decomposition theorem (Proposition 5.4 of \cite{Ca}), which provides a decomposition of the universal orbifold covering $\widetilde{X}$ of $X$ as a product $M_{0}\times M_{1}\times\cdots\times M_{k}$, where $M_{0}$ is a Euclidean space, and $M_{1},\cdots,M_{k}$ have all irreducible holonomy representation. The orbifold version of the Cheeger-Gromoll theorem is provided by Borzellino and Zhu (see \cite{BZ}), and implies that $M_{1},\cdots,M_{k}$ are all compact and their Ricci curvature is zero. The remaining part of the proof is similar to the one for compact manifolds.

\subsubsection{Related notions and examples}

There are other definitions which are related to irreducible symplectic orbifolds. 

\begin{definition}
A \textit{primitively symplectic V-manifold} is a symplectic orbifold $X$ such that $h^{2,0}(X)=1$. An \textit{irreducible symplectic V-manifold} is a compact, connected K\"ahler orbifold $X$ which is simply connected and such that $h^{2,0}(X)=1$.
\end{definition}

Primitively symplectic V-manifolds are introduced in \cite{Fu}, section 2.1, while irreducible symplectic V-manifolds are introduced in \cite{MT}, Definition 1.3.(iv). Clearly an irreducible symplectic V-manifold is a compact, connected K\"ahler primitively symplectic V-manifold which is simply connected. Moreover, an irreducible symplectic orbifold is an irreducible symplectic V-manifold. 

We will discuss in the following several examples of irreducible symplectic V-manifolds/orbifolds of dimension 4 and 6 which appear in the literature. In particular, we will see that irreducible symplectic V-manifolds are not always irreducible symplectic orbifolds.\par\bigskip

\paragraph{\textbf{Symmetric products of K3 surfaces}.}
\hfill\\
\hfill\\
If $S$ is a K3 surface and $m\geq 2$, the symmetric product $X:=Sym^{m}(S)$ is an irreducible symplectic V-manifold. Indeed, it is a compact, connected K\"ahler orbifold having $Y:=Hilb^{m}(S)$ as a resolution of the singularities. As $Y$ is an irreducible symplectic manifold, it is simply connected and has $h^{2,0}(Y)=1$: it follows that $X$ is simply connected and $h^{2,0}(X)=1$ (see as instance Proposition 2.13 of \cite{Fu}). 

We notice that the smooth locus $X^{s}$ of $X$ is not simply connected, since $\pi_{1}(X^{s})\simeq \Sigma_{m}$. It follows that $X$ is not an irreducible symplectic orbifold.\par\bigskip

\paragraph{\textbf{Fujiki's examples}.}
\hfill\\
\hfill\\
Tables I and II in section 13 of \cite{Fu} present a list of 18 examples of irreducible symplectic V-manifolds of dimension 4 constructed as follows: let $S$ be a K3 surface or a $2-$dimensional complex torus, $H$ a finite group acting symplectically on $S$ (i. e. each element $h\in H$ preserves the holomorphic symplectic form of $S$) and $\tau$ an automorphism of $H$ of order 2. Then $H$ acts on $S\times S$ via the action mapping $(h,(s,t))$ to $(h(s),\tau(h)(t))$. Let $G(H)$ be the subgroup of $Aut(S\times S)$ generated by $H$ and by the involution $\iota$ mapping $(s,t)$ to $(t,s)$, and consider $Y:=S\times S/G(H)$.

The fixed locus of the action of $G(H)$ may have $2-$dimensional components, and apart from that it has only isolated points: the singular locus of $Y$ is then given by a (possibly empty) $2-$dimensional locus $\Sigma$ and by a finite number of points. Blowing up $\Sigma$, one gets a V-manifold $X$ whose singular locus is given by a finite number of points. Theorem 13.1 of \cite{Fu} shows that $X$ is an irreducible symplectic V-manifold.

We reproduce here Tables I and II for the convenience of the reader: here $b_{2}$ is the second Betti number of the variety, and $a_{k}$ is the number of singular points of type $\widehat{A}_{k}=(\mathbb{C}^{4}/g_{k},0)$, where $g_{k}=(\zeta_{k},\zeta_{k},\zeta_{k}^{-1},\zeta_{k}^{-1})$ and $\zeta_{k}:=e^{\frac{2\pi i}{k}}$ (the sum $a_{2}+a_{3}+a_{4}+a_{6}$ is the number of singular points of the variety). The varieties $X_{p}$ are obtained starting from a K3 surface $S$, while the varieties $Y_{n}$ are obtained starting from a $2-$dimensional complex torus.
\begin{table}
\label{table:Table1}
\caption{Fujiki's examples}
\begin{tabular}{|c|c|c|c|c|c|c|}
\hline
Symbol & $H$ & $b_{2}$ & $a_{2}$ & $a_{3}$ & $a_{4}$ & $a_{6}$ \\
\hline
\hline
$X_{1}$ & $id_{S}$ & 23 & 0 & 0 & 0 & 0 \\
\hline
$X_{2}$ & $\mathbb{Z}/2\mathbb{Z}$ & 16 & 28 & 0 & 0 & 0 \\
\hline
$X_{3}$ & $(\mathbb{Z}/2\mathbb{Z})^{\oplus 2}$ & 14 & 36 & 0 & 0 & 0 \\
\hline
$X_{4}$ & $(\mathbb{Z}/2\mathbb{Z})^{\oplus 3}$ & 16 & 28 & 0 & 0 & 0 \\
\hline
$X_{5}$ & $\mathbb{Z}/3\mathbb{Z}$ & 11 & 0 & 15 & 0 & 0 \\
\hline
$X_{6}$ & $(\mathbb{Z}/2\mathbb{Z})^{\oplus 2}$ & 7 & 0 & 12 & 0 & 0 \\
\hline
$X_{7}$ & $\mathbb{Z}/4\mathbb{Z}$ & 10 & 10 & 0 & 6 & 0 \\
\hline
$X_{8}$ & $(\mathbb{Z}/4\mathbb{Z})^{\oplus 2}$ & 8 & 12 & 0 & 0 & 0 \\
\hline
$X_{9}$ & $\mathbb{Z}/6\mathbb{Z}$ & 8 & 7 & 6 & 0 & 1 \\
\hline
$X_{10}$ & $(\mathbb{Z}/2\mathbb{Z})\oplus(\mathbb{Z}/4\mathbb{Z})$ & 10 & 16 & 0 & 4 & 0 \\
\hline
$X_{11}$ & $(\mathbb{Z}/2\mathbb{Z})\oplus(\mathbb{Z}/6\mathbb{Z})$ & 8 & 6 & 6 & 0 & 0 \\ 
\hline
$Y_{1}$ & $\mathbb{Z}/3\mathbb{Z}$ & 7 & 0 & 36 & 0 & 0 \\
\hline
$Y_{2}$ & $(\mathbb{Z}/3\mathbb{Z})^{\oplus 2}$ & 7 & 0 & 27 & 0 & 0 \\
\hline
$Y_{3}$ & $(\mathbb{Z}/3\mathbb{Z})^{\oplus 3}$ & 7 & 0 & 0 & 0 & 0 \\
\hline
$Y_{4}$ & $\mathbb{Z}/4\mathbb{Z}$ & 8 & 54 & 0 & 6 & 0 \\
\hline
$Y_{5}$ & $(\mathbb{Z}/4\mathbb{Z})\oplus(\mathbb{Z}/2\mathbb{Z})$ & 10 & 52 & 0 & 4 & 0 \\
\hline
$Y_{6}$ & $(\mathbb{Z}/4\mathbb{Z})\oplus(\mathbb{Z}/2\mathbb{Z})^{\oplus 2}$ & 14 & 48 & 0 & 0 & 0 \\
\hline
$Y_{7}$ & $\mathbb{Z}/6\mathbb{Z}$ & 8 & 35 & 16 & 0 & 0 \\
\hline
\end{tabular}
\end{table}

In all the previous examples, the singular locus of $S\times S/G(H)$ has always at least one irreducible component of dimension 2: we have 1 component for $X_{1}$, $X_{5}$, $X_{6}$, $Y_{1}$, $Y_{2}$ and $Y_{3}$; 2 components for $X_{2}$, $X_{7}$, $X_{9}$, $Y_{4}$ and $Y_{7}$; 4 components for $X_{3}$, $X_{8}$, $X_{10}$, $X_{11}$ and $Y_{5}$; 8 components for $X_{4}$ and $Y_{6}$.

The varieties $X_{1}$ and $Y_{3}$ are the only smooth examples we get, and they are both irreducible symplectic manifolds: $X_{1}$ is the Hilbert scheme of two points on $S$, while $Y_{3}$ is deformation equivalent to $Kum^{2}(T)$ (see Remark 13.2.3 and Proposition 14.3 of \cite{Fu}). 

If $p\neq 1$ and $n\neq 3$, then $X_{p}$ and $Y_{n}$ are singular symplectic varieties whose singular locus has codimension 4: by Corollary 1 of \cite{N1} they all have terminal singularities, and by the Main Theorem of \cite{N2} their deformations are locally trivial. It follows that if $p\neq q$ and $(p,q)\neq(2,4)$ then $X_{p}$ and $X_{q}$ are not deformation equivalent; if $n\neq m$ then $Y_{n}$ is not deformation equivalent to $Y_{m}$, and for every $p,n$ we have that $X_{p}$ is not deformation equivalent to $Y_{n}$. 

It is not known if $X_{2}$ and $X_{4}$ give different deformation classes. Anyway, the Fujiki examples provide at least 17 different deformation classes of irreducible symplectic V-manifolds in dimension 4, 15 of which are singular. The values of $b_{2}$ of these examples are 7, 8, 10, 11, 14, 16 and 23. It is still an open problem to determine if $X_{p}$ and $Y_{n}$ are irreducible symplectic orbifolds if $p\neq 1$ and $n\neq 3$.\par\bigskip

\paragraph{\textbf{Quotients of $Hilb^{n}(K3)$}.}
\hfill\\
\hfill\\
We now consider quotients of $Hilb^{n}(S)$ where $S$ is a projective K3 surface. We start with the case $n=2$, so that the fixed locus $Fix(G)$ of $G$ is a union of a finite number of points (depending on the order of $G$) and a (possibly empty) holomorphic symplectic surface, and the same holds for the singular locus of $Hilb^{2}(S)/G$. By blowing up the singular surface one then gets an irreducible symplectic V-manifold $M_{G}$ (see \cite{Fu}).

A particular case is when $G=\langle\phi\rangle$, where $\phi$ is a symplectic automorphism of prime order. By Corollary 2.13 of \cite{Mon2}, the order $p$ of $\phi$ can only be 2, 3, 5, 7 or 11. The case of $p=2$, i. e. $\phi$ is a symplectic involution, is studied in \cite{MT} and \cite{Men2}. The case $p=3$ is studied in \cite{Men2}. The case $p=11$ is studied in \cite{Men1}. The cases $p=5,7$ will be treated by Menet in a forthcoming paper.

\begin{example}
{\rm If $G=\langle\phi\rangle$ where $\phi$ is a symplectic involution, then $M_{2}:=M_{G}$ is an irreducible symplectic V-manifold with 28 isolated singular points, and by Theorem 2.5 of \cite{Men2} we have $b_{2}(M_{2})=16$. By Corollary 3.23 of \cite{Men3} and Proposition 5.1 of \cite{MT} we have that the topological Euler number of $M_{2}$ is 268.}

{\rm It follows that $M_{2}$ is not deformation equivalent to any of the Fujiki examples above: the only examples in Table 1 having 28 singular points are $X_{2}$ and $X_{4}$ (for both of which $b_{2}=16$), whose topological Euler number is 226 (see Remark 13.2.4 of \cite{Fu}).}
\end{example}

\begin{example}
{\rm If $G=\langle\phi\rangle$ where $\phi$ is a natural symplectic automorphism of order 3, then $Fix(G)$ is given by 27 isolated points. We get that $M_{3}:=M_{G}$ is an irreducible symplectic V-manifold with 27 isolated singular points.}

{\rm This example is then not deformation equivalent to $M_{2}$, and neither to any of the Fujiki examples in Table 1: the only one of that list having 27 singular points is $Y_{2}$, whose $b_{2}=7$. But by Theorem 1.3 of \cite{Men1} we have $b_{2}(M_{3})=11$, so $M_{3}$ is a new deformation class.}
\end{example}

\begin{example}
{\rm Examples 4.5.1 and 4.5.2 of \cite{Mon2} provide two K3 surfaces $S_{1}$ and $S_{2}$ such that $Hilb^{2}(S_{i})$ has an automorphism $\sigma_{i}$ of order 11. Let $G_{i}:=\langle\sigma_{i}\rangle$: in both cases $Fix(G_{i})$ is given by 5 isolated points. Then $M_{11}^{i}:=M_{G_{i}}$ is an irreducible symplectic V-manifold with 5 singular points.}

{\rm By \cite{Men1} we have $b_{2}(M_{11}^{i})=3$ and that $M^{1}_{11}$ and $M^{2}_{11}$ are not deformation equivalent, so they provide two more deformation classes.}
\end{example}

We then get 4 more deformation classes of singular irreducible symplectic V-manifolds of dimension 4 to add to the 15 classes presented by Fujiki. If we now consider a K3 surface $S$ and $Hilb^{n}(S)$ for $n\geq 3$, let $\phi$ be a symplectic involution on $Hilb^{n}(S)$. The quotient $H_{n}$ of $Hilb^{n}(S)$ by the action of $\phi$ is again an irreducible symplectic $V-$manifold (see \cite{Fu}, or Proposition \ref{proposition:quotisv} below).

We are in the position to describe which of these examples are irreducible symplectic orbifolds, as the following shows:

\begin{proposition}
\label{proposition:quothilb}
Let $S$ be a projective K3 surface and $G$ a finite group of automorphisms of $Hilb^{n}(S)$ acting symplectically.
\begin{enumerate}
 \item If the codimension of $Fix(G)$ is at least 4, then $M_{G}:=Hilb^{2}(S)/G$ is not an irreducible symplectic orbifold.
 \item If $n=2$ and $G$ is generated by a symplectic involution, then the singular locus of $Hilb^{2}(S)/G$ is given by 28 isolated points and a K3 surface $\Sigma$. Then the partial resolution $M_{G}$ of $Hilb^{2}(S)/G$ obtained by blowing-up $\Sigma$ is an irreducible symplectic orbifold.
\end{enumerate}
In particular $M_{2}$ is an irreducible symplectic orbifold, while $M_{3}$, $M_{11}^{1}$, $M_{11}^{2}$ and $H_{n}$ are irreducible symplectic $V-$manifolds which are not irreducible symplectic orbifolds.
\end{proposition}

\proof Suppose first that $Fix(G)$ has codimension at least 4 in $Hilb^{n}(S)$, so the singular locus of $M_{G}$ has codimension at least 4. We know that $M_{G}$ is an irreducible symplectic V-manifold: we let $U$ be its smooth locus and $f:Hilb^{n}(S)\longrightarrow M_{G}$ the quotient morphism. Hence $f:f^{-1}(U)\longrightarrow U$ is a finite \'etale covering whose degree is the order of $G$: hence $U$ is not simply connected, and $M_{G}$ is not an irreducible symplectic orbifold.

We notice that for $M_{3}$, $M_{11}^{1}$ and $M_{11}^{2}$ we have that $Fix(G)$ is given by isolated points (see Example 2 and Example 3 above), so they are not irreducible symplectic orbifolds. The fixed locus $Fix(G)$ in the case of $H_{n}$ has codimension at least 4 by Theorem 1.1 of \cite{KMO}, so that $H_{n}$ is not an irreducible symplectic orbifold.

The case of $n=2$ and $G$ generated by a symplectic involution is due to Menet (see Remark 3.22 of \cite{Men3}). We need to show that the smooth locus $U$ of $M_{G}$ is simply connected. Let $D$ be the exceptional divisor of $M_{G}$ coming from the blow-up $b:M_{G}\longrightarrow X_{G}:=Hilb^{2}(S)/G$.

Now, let $\Sigma'$ be the $2-$dimensional component of the singular locus of $X_{G}$, and $U'$ the smooth locus of $X_{G}$. We notice that the blow-up morphism $b$ gives an isomorphism between $U\setminus D$ and $U'$. Moreover the quotient morphism $f:Hilb^{2}(S)\longrightarrow X_{G}$ gives a double \'etale covering $Hilb^{2}(S)\setminus Fix(G)\longrightarrow U'$. As $Fix(G)$ has codimension 2 in $Hilb^{2}(S)$, and $Hilb^{2}(S)$ is smooth, it follows that $\pi_{1}(Hilb^{2}(S))\simeq\pi_{1}(Hilb^{2}(S)\setminus Fix(G))$, so that $Hilb^{2}(S)\setminus Fix(G)$ is simply connected.

As a consequence, we see that $\pi_{1}(U\setminus D)\simeq\pi_{1}(U')\simeq\mathbb{Z}/2\mathbb{Z}$. Hence, there is an étale covering $\pi':Y'\longrightarrow U\setminus D$ of degree 2, which extends to a finite covering $\pi:Y\longrightarrow U$ branched along $D$. Notice that $Y'$ is simply connected, hence it follows that $Y$ is simply connected as well.

Let now $x_{0}\in D$ and consider a loop $\gamma$ in $U$ pointed at $x_{0}$. Let $\widetilde{\gamma}$ be a lift of $\gamma$ to $Y$: notice that as $x_{0}$ is a branching point, the fiber of $\pi$ over $x_{0}$ is given by a unique point $y_{0}$, so $\widetilde{\gamma}$ is a loop in $Y$ pointed at $y_{0}$. Since $Y$ is simply connected, it follows that the homotopy class of $\widetilde{\gamma}$ is zero in $\pi_{1}(Y,y_{0})$, so the homotopy class of $\gamma$ is zero as well in $\pi_{1}(U,x_{0})$, and we are done.\endproof

\paragraph{\textbf{Quotients of $Kum^{2}(S)$}.}
\hfill\\
\hfill\\
Similar considerations are done in \cite{KaMe} for quotients of $Kum^{2}(S)$ for an Abelian surface $S$. If $\sigma$ is a symplectic involution on $Kum^{2}(S)$, by Theorem 7.5 of \cite{KaMe} the fixed locus of $\sigma$ is given by 36 points together with a K3 surface $\Sigma$. 

The quotient $Kum^{2}(S)/\sigma$ has then a singular locus given by 36 points and a K3 surface. Blowing up the image of $\Sigma$ in $Kum^{2}(S)/\sigma$ one gets a $V-$manifold $K$ which has 36 singular points. Again $K$ is an irreducible symplectic V-manifold, and $b_{2}(K)=8$ (see \cite{KaMe}), hence $K$ is not deformation equivalent to $M_{2}$, $M_{3}$, $M_{11}^{1}$ and $M_{11}^{2}$. 

Moreover, the only examples in Table 1 having 36 singular points are $X_{3}$ and $Y_{1}$, whose $b_{2}$ is 14 and 7 respectively, so $K$ is not deformation equivalent to any of the Fujiki examples: we then have another deformation class of singular irreducible symplectic V-manifolds to add the the previous 19 classes. The same proof of point 2 of Proposition \ref{proposition:quothilb} gives the following:

\begin{proposition}
The variety $K$ is an irreducible symplectic orbifold.
\end{proposition}

\paragraph{\textbf{Examples of Markushevich-Tikhomirov}.}
\hfill\\
\hfill\\
Markushevich and Tikhomirov present in \cite{MT} a different construction of irreducible symplectic V-manifolds of dimension 4. Let $X$ be a del Pezzo surface which is the double cover of $\mathbb{P}^{2}$ branched along a smooth quartic $B_{0}$ with 28 bitangent lines, and let $S$ be the double cover of $X$ branched along the curve $\Delta_{0}$ such that $\Delta_{0}+i(\Delta_{0})$ is the inverse image of a smooth quartic curve of $\mathbb{P}^{2}$ which is totally tangent to $B_{0}$ at eight distinct points (here $i$ is the involution on $X$ induced by the double cover $X\longrightarrow\mathbb{P}^{2}$). 

The surface $S$ is a K3 surface, and let $M_{k}$ be the moduli space of torsion sheaves on $S$ with first Chern class $H$ (the pull-back of $-K_{X}$) and Euler character $k-2$ for $k\in 2\mathbb{Z}$, i. e. of Mukai vector $v=(0,H,k-2)$. We notice that $M_{k}$ is a projective variety of dimension 6: as $k$ is even, $M_{k}$ has exactly 28 singular points (see Proposition 1.12.$(iii)$ of \cite{MT}). Moreover, the moduli space $M_{k}$ has an involution $\sigma$ mapping $L\in M_{k}$ to $\mathcal{E}xt^{1}_{\mathcal{O}_{S}}(L,\mathcal{O}_{S}(-H))$. 

If $\tau$ is the involution on $S$ induced by the double cover $S\longrightarrow X$, we let $\kappa:=\tau^{*}\circ\sigma$: the fixed locus of $\kappa$ has a $4-$dimensional irreducible component, denoted $\mathcal{P}^{k}$. The morphism mapping $L\in\mathcal{P}^{k}$ to $L(H)\in \mathcal{P}^{k+2}$ is an isomorphism, hence we get at most two non-isomorphic varieties $\mathcal{P}^{0}$ and $\mathcal{P}^{2}$. 

By \cite{MT} we know that $\mathcal{P}^{0}$ and $\mathcal{P}^{2}$ are both irreducible symplectic V-manifolds having 28 singular points. By Lemma 5.2 and Corollary 5.7 of \cite{MT} we have that $\mathcal{P}^{0}$ is birational to $M_{2}$ via a Mukai flop, and by Corollary 3.23 of \cite{Men3}: it follows that $\mathcal{P}^{0}$ is an irreducible symplectic orbifold. 

It is not known if $\mathcal{P}^{0}$ and $\mathcal{P}^{2}$ are birational, isomorphic nor deformation equivalent, and it is still an open question if $\mathcal{P}^{2}$ is an irreducible symplectic orbifold.\par\bigskip

\paragraph{\textbf{Example of Matteini}.}
\hfill\\
\hfill\\
A similar construction to that of Markushevich and Tikhomirov is presented in \cite{Mat}. Take a K3 surface $S$ which is a double cover of a generic cubic surface $Y$ with involution $\tau$, and take $M$ to be the moduli space of semistable torsion sheaves with first Chern class $H$ (the pull-back of $-K_{Y}$) and Euler character $-3$: this is a singular projective variety of dimension 8.

One still has an involution $\sigma$ on $M$ obtained as before, and we let $\kappa:=\tau^{*}\circ\sigma$. The fixed locus of $\kappa$ has a $6-$dimensional irreducible component $\mathcal{P}$, whose singular locus is the union of 27 singular K3 surfaces. In \cite{Mat} it is shown that $\mathcal{P}$ is an irreducible symplectic V-manifold of dimension 6. It is not known if this example is an irreducible symplectic orbifold.

\subsection{Varieties with canonical singularities}

A further generalization of the Beauville-Bogomolov decomposition theorem was obtained for projective varieties with canonical singularities by the works of Druel, Greb, Guenancia, H\"oring, Kebekus and Peternell, in particular \cite{GKP}, \cite{D}, \cite{DG}, \cite{GGK} and \cite{HP}. 

\subsubsection{The decomposition for singular projective varieties}

We introduce the following notation: if $X$ is a normal variety and $X_{reg}$ is the smooth locus of $X$ whose open embedding in $X$ is $j:X_{reg}\longrightarrow X$, for every $p\in\mathbb{N}$ such that $0\leq p\leq\dim(X)$ we let $$\Omega_{X}^{[p]}:=j_{*}\Omega^{p}_{X_{reg}}=\big(\wedge^{p}\Omega_{X}\big)^{**},$$whose global sections are called \textit{reflexive $p-$forms} on $X$. 

Notice that by definition of $\Omega_{X}^{[p]}$ we have $H^{0}(X,\Omega_{X}^{[p]})=H^{0}(X_{reg},\Omega_{X_{reg}}^{p})$. Theorem 1.4 of \cite{GKKP} shows that if $X$ is a quasi-projective variety with klt singularities and $\pi:\widetilde{X}\longrightarrow X$ is a log-resolution, then for every $p\in\mathbb{N}$ such that $0\leq p\leq\dim(X)$ the sheaf $\pi_{*}\Omega^{p}_{\widetilde{X}}$ is reflexive. This implies in particular that $H^{0}(X,\Omega_{X}^{[p]})\simeq H^{0}(\widetilde{X},\Omega^{p}_{\widetilde{X}})$ (see Observation 1.3 therein).

As shown in \cite{GKKP}, if $f:Y\longrightarrow X$ is a finite, dominant morphism between two irreducible normal varieties, then there is a morphism of reflexive sheaves $f^{*}\Omega_{X}^{[p]}\longrightarrow\Omega_{Y}^{[p]}$, induced by the usual pull-back morphism of forms on the smooth loci, giving a morphism $f^{[*]}:H^{0}(X,\Omega_{X}^{[p]})\longrightarrow H^{0}(Y,\Omega_{Y}^{[p]})$, called \textit{reflexive pull-back morphism}.

We recall the definitions of symplectic form and symplectic variety (see \cite{B2}).

\begin{definition}
Let $X$ be a normal variety.
\begin{enumerate}
 \item A \textit{symplectic form} on $X$ is a closed reflexive $2-$form $\sigma$ on $X$ which is non-degenerate at each point of $X_{reg}$.
 \item If $\sigma$ is a symplectic form on $X$, the pair $(X,\sigma)$ is a \textit{symplectic variety} if for every resolution $f:\widetilde{X}\longrightarrow X$ of the singularities of $X$, the holomorphic symplectic form $\sigma_{reg}:=\sigma_{|X_{reg}}$ extends to a holomorphic $2-$form on $\widetilde{X}$.
 \item If $(X,\sigma)$ is a symplectic variety and $f:\widetilde{X}\longrightarrow X$ is a resolution of the singularities over which $\sigma_{reg}$ extends to a holomorphic symplectic form, we say that $f$ is a \textit{symplectic resolution}.
\end{enumerate}
\end{definition}

A normal variety having a symplectic form and whose singular locus has codimension at least 4 is a symplectic variety (see \cite{Fle}), and a symplectic variety has terminal singularities if and only if its singular locus has codimension at least 4 (Corollary 1 of \cite{N1}).

We now define irreducible Calabi-Yau and irreducible symplectic varieties following \cite{GKP}. If $X$ and $Y$ are two irreducible normal projective varieties, a \textit{finite quasi-\'etale morphism} $f:Y\longrightarrow X$ is a finite morphism which is \'etale in codimension one.

\begin{definition}
\label{defn:irrvar}
Let $X$ be an irreducible normal projective variety with trivial canonical divisor and canonical singularities, of dimension $d\geq 2$.
\begin{enumerate}
 \item The variety $X$ is \textit{irreducible Calabi-Yau} if for every $0<p<d$ and for every finite quasi-\'etale morphism $Y\longrightarrow X$, we have $H^{0}(Y,\Omega_{Y}^{[p]})=0$.
 \item The variety $X$ is \textit{irreducible symplectic} if it has a symplectic form $\sigma\in H^{0}(X,\Omega_{X}^{[2]})$, and for every finite quasi-\'etale morphism $f:Y\longrightarrow X$ the exterior algebra of reflexive forms on $Y$ is spanned by $f^{[*]}\sigma$.
\end{enumerate}
\end{definition}

A description of irreducible Calabi-Yau and irreducible symplectic varieties in terms of holonomy is available by the work of Greb, Guenancia and Kebekus \cite{GGK}. More precisely, if $H$ is an ample divisor on $X$, by \cite{EGZ} there is a singular Ricci-flat K\"ahler metric $\omega_{H}$ in $c_{1}(H)$, inducing a Riemannian metric $g_{H}$ on $X_{reg}$. 

We let $Hol(X_{reg},g_{H})$ be the holonomy group of this metric. Proposition F of \cite{GKP} shows that a normal projective variety of dimension $n$ with klt singularities and trivial canonical bundle is an irreducible Calabi-Yau variety if and only if $Hol(X_{reg},g_{H})$ is isomorphic to $SU(n)$, and it is an irreducible symplectic variety if and only if $Hol(X_{reg},g_{H})$ is isomorphic to $Sp(n/2)$.

The decomposition theorem for singular projective varieties is the following:

\begin{theorem}
\label{theorem:bbdecompsing}
Let $X$ be a normal projective varieties with klt singularities and numerically trivial canonical bundle. Then $X$ has a finite quasi-\'etale covering $f:Y\longrightarrow X$, where $Y$ is a normal projective variety with canonical singularities which is a product of complex tori, irreducible Calabi-Yau varieties and irreducible symplectic varieties.
\end{theorem}

This is Theorem 1.15 of \cite{HP}, and the proof can be found therein. It consists of three major parts: one is the holonomy decomposition obtained by Greb, Guenancia and Kebekus in \cite{GGK} (namely Theorem B and Proposition D therein); a second one is an algebraic integrability theorem of Druel, which is Theorem 1.4 of \cite{D}; the final ingredient is Theorem 1.1 of \cite{HP}. Less general versions of the Bogomolov decomposition theorem in the projective singular setting were previously obtained in \cite{GKP}, \cite{D} and \cite{DG}. The complete proof can be found in section 4 of \cite{HP} (the proof of Theorem 1.5 therein).

\begin{remark}
{\rm We notice that even if the statement of Theorem 1.15 of \cite{HP} gives the existence of a quasi-\'etale covering $f:Y\longrightarrow X$, this is consistent with the statement of Theorem \ref{theorem:bbdecompsing} since they define quasi-\'etale morphisms as finite morphisms whose ramification divisor is empty.}
\end{remark}

\subsubsection{Relation between the previous notions}

The first result we state is about the relation between irreducible symplectic manifolds, orbifolds and varieties.

\begin{proposition}
\label{proposition:ismov}
The following properties hold.
\begin{enumerate}
 \item Irreducible symplectic manifolds are irreducible symplectic orbifolds.
 \item Projective irreducible symplectic orbifolds are irreducible symplectic varieties.
 \item Smooth irreducible symplectic varieties are irreducible symplectic manifold.
\end{enumerate}
\end{proposition}

\proof The first point is trivial, and the last is a consequence of Proposition A.1 of \cite{HNW}.

Suppose that $X$ is a projective irreducible symplectic orbifold. Then $X$ is an normal projective variety with rational Cohen-Macaulay singularities and trivial canonical bundle. It follows that $X$ has rational Gorenstein singularities. Moreover, it has a holomorphic symplectic form on its singular locus, so by Theorem 6 of \cite{N3} it follows that $X$ is a projective symplectic variety. In particular it has canonical singularities (see \cite{B2}).

By Theorem \ref{theorem:bbdecompsing} there is then a finite quasi-\'etale covering $f:Y\longrightarrow X$, where $Y$ is a product of Abelian varieties, irreducible Calabi-Yau varieties and irreducible symplectic varieties. As $X_{reg}$ is simply connected by definition of irreducible symplectic orbifold, it follows that $f$ is an isomorphism. As a consequence $Y$ is simply connected, so it has no factor which is an Abelian variety, and it is a symplectic variety, hence it has no factor which is an irreducible Calabi-Yau variety. 

Hence $Y$ is a product of irreducible symplectic varieties. But as $H^{0}(X,\Omega_{X}^{[2]})$ is one dimensional by definition of irreducible symplectic orbifold, the same holds for $Y$. If $Y$ is a product of $m$ irreducible symplectic varieties, we have that $H^{0}(Y,\Omega_{Y}^{[2]})$ has dimension $m$: it follows that $m=1$, so $Y$ is an irreducible symplectic variety. As $X$ is isomorphic to $Y$, we are done.\endproof

Irreducible symplectic V-manifolds are not necessarily irreducible symplectic varieties. A first example of this is given by symmetric products of K3 surfaces: if $S$ is a K3 surface and $m\in\mathbb{N}$, $m\geq 2$, then $X:=Sym^{m}(S)$ is an irreducible symplectic V-manifold. Anyway it has a finite quasi-\'etale covering $S^{m}\longrightarrow X$, and $h^{0}(S^{m},\Omega^{2}_{S^{m}})=m$, so $X$ is not an irreducible symplectic variety.

By Theorem I of \cite{GKP} we know that all irreducible symplectic varieties are simply connected, so if $X$ is a primitively symplectic V-manifold which is an irreducible symplectic variety, then $X$ is an irreducible symplectic V-manifold. Anyway, as symplectic singularities are not, in general, quotient singularities, we cannot expect that an irreducible symplectic variety is an irreducible symplectic V-manifold: an example will be given in the last section.

The following is a criterion to guarantee that some quotients of an irreducible symplectic manifold are irreducible symplectic varieties (and irreducible symplectic $V-$manifolds as well).

\begin{proposition}
\label{proposition:quotisv}Let $X$ be an irreducible symplectic manifold, $G\subseteq Aut(X)$ a finite subgroup acting symplectically on $X$ and $Y:=X/G$.
\begin{enumerate}
 \item The quotient $Y$ is an irreducible symplectic $V-$manifold.
 \item If $Y$ has terminal singularities and $X$ is projective, then $Y$ is an irreducible symplectic variety.
\end{enumerate}
\end{proposition}

\proof The proof of the first part is basically contained in \cite{Fu}, but we present it here for the reader's sake. The fact that $X$ is a compact, connected K\"ahler orbifold is Lemma 1.4 of \cite{Fu}, and by Lemma 2.4 of \cite{Fu} we have that $Y$ is symplectic. As $X$ is simply connected, by Lemma 1.2 of \cite{Fu} it follows that $Y$ is simply connected as well. Finally, we have $H^{2,0}(Y)\simeq H^{2,0}(X)^{G}$, where $H^{2,0}(X)^{G}$ is the space of $G-$invariant sections. But as $G$ acts symplectically and $X$ is an irreducible symplectic manifold, it follows that $H^{2,0}(Y)$ is $1-$dimensional, so that $Y$ is an irreducible symplectic $V-$manifold.

For the second part, as we know that $Y$ is a normal projective symplectic variety with terminal singularities, in order to show that it is an irreducible symplectic variety we just need to look at its finite quasi-\'etale coverings. So, let $f:Y'\longrightarrow Y$ be a finite quasi-\'etale covering of $Y$. As the singularities of $Y$ are terminal and as $X$ is simply connected (being an irreducible symplectic manifold), it follows that $X$ is the universal covering of $Y$.

In particular, it follows that we have a finite quasi-\'etale covering $\pi':X\longrightarrow Y'$ given by the quotient by a subgroup $G'$ of $G$, and if we let $\pi:X\longrightarrow Y$, then $\pi=f\circ\pi'$. As $X$ is an irreducible symplectic manifold, its exterior algebra of holomorphic forms is spanned by the symplectic form on $X$, which is the reflexive pull-back of the one on $Y$. 

Moreover, we notice that $H^{0}(Y',\Omega_{Y'}^{[p]})\simeq H^{p,0}(X)^{G'}$, and as $X$ is an irreducible symplectic manifold and $G$ (and hence $G'$) acts symplectically we see that if $p$ is even then $H^{0}(Y',\Omega_{Y'}^{[p]})$ is spanned by $f^{[*]}\sigma^{p/2}$ (where $\sigma$ is the symplectic form on $Y$) and if $p$ is odd then $H^{0}(Y',\Omega_{Y'}^{[p]})=0$. But this shows that $Y$ is an irreducible symplectic variety.\endproof 

As an application of this, we see that the examples $M_{3}$, $M_{11}^{1}$, $M_{11}^{2}$ and $H_{n}$ presented in section 2.1.2 are all examples of irreducible symplectic $V-$manifolds which are irreducible symplectic varieties (but not irreducible symplectic orbifolds). 

The partial resolution $M_{2}$ of the quotient of $Hilb^{2}(S)$ by a symplectic involution (where $S$ is a K3 surface), the partial resolution $K_{2}$ of the quotient of $Kum^{2}(T)$ by a symplectic involution (where $T$ is a $2-$dimensional complex torus), and the example $\mathcal{P}^{0}$ of \cite{MT} (which is deformation equivalent to $M_{2}$) are all irreducible symplectic varieties as they are irreducible symplectic orbifolds.  

For all other examples in section 2.1.2 (those of Fujiki, the example $\mathcal{P}^{2}$ of \cite{MT} and the example of Matteini), it is not known if they are irreducible symplectic varieties.

\subsubsection{Related notions}

Irreducible symplectic varieties appear in several papers under different definitions. A first one appears in \cite{BL}, and it is defined as follows (the name given to these varieties in \cite{BL} is \textit{irreducible symplectic varieties}).

\begin{definition}
A \textit{resolvable symplectic variety} is a normal, compact K\"ahler space $X$ whose smooth locus has a holomorphic symplectic form $\sigma$, and which has a symplectic resolution of the singularities which is an irreducible symplectic manifold.
\end{definition} 

A projective resolvable symplectic variety is not always an irreducible symplectic variety: if $S$ is a projective K3 surface and $m\geq 2$, then $Sym^{m}(S)$ is a projective resolvable symplectic variety (since it is a normal projective symplectic variety having $Hilb^{m}(S)$ as a symplectic resolution), but it is not an irreducible symplectic variety.

Similarly, singular Kummer surfaces (i. e. a surface $S$ obtained as quotient of an Abelian surface $A$ by the involution mapping $p\in A$ to $-p\in A$) are resolvable symplectic surfaces (they have a symplectic resolution which is a K3 surface), irreducible symplectic V-manifolds but not irreducible symplectic varieties (as the quotient map $A\longrightarrow S$ is a finite quasi-\'etale covering and $h^{1,0}(A)\neq 0$).

We will see in section 2.2.4 examples of irreducible symplectic varieties which are not resolvable symplectic varieties, and of resolvable symplectic varieties which are not irreducible symplectic V-manifolds. Anyway we have the following:

\begin{proposition}
\label{proposition:symplbl}
If $X$ is an irreducible symplectic variety (resp. an irreducible symplectic V-manifold) having a symplectic resolution $Y$, then $Y$ is an irreducible symplectic manifold. In particular, $X$ is a resolvable symplectic variety.
\end{proposition}

\proof If $X$ is an irreducible symplectic variety, this is Remark 1.16 of \cite{PR3}. If $X$ is an irreducible symplectic V-manifold, then $X$ has canonical singularities. By \cite{T} we get $\pi_{1}(X)\simeq\pi_{1}(Y)$, so $Y$ is simply connected. Finally, by Theorem 1.4 of \cite{GKKP} we have $h^{0}(Y,\Omega_{Y}^{2})=h^{0}(X,\Omega_{X}^{[2]})$, which is 1 by definition, and we are done.\endproof

A further definition of irreducible symplectic variety appears in \cite{S}, where it is defined as a projective symplectic variety $X$ such that $h^{1}(X,\mathcal{O}_{X})=0$ and $h^{0}(X,\Omega_{X}^{[2]})=1$. It is called Namikawa symplectic variety if it is moreover $\mathbb{Q}-$factorial and its singular locus has codimension at least 4 (see Definition 1 therein).

\begin{definition}
We will call \textit{Namikawa symplectic variety} a normal, compact K\"ahler complex space $X$ such that $h^{1}(X,\mathcal{O}_{X})=0$ and $h^{0}(X,\Omega_{X}^{[2]})=1$.
\end{definition}

Namikawa symplectic varieties are the most general kind of varieties we will deal with. Namely:

\begin{proposition}
\label{proposition:nami}
Irreducible symplectic varieties, resolvable symplectic varieties and irreducible symplectic V-manifolds are all Namikawa symplectic varieties.
\end{proposition}

\proof The proof that irreducible symplectic and resolvable symplectic varieties are Namikawa symplectic is given in Propositions 1.9 and 1.10 of \cite{PR3}. If $X$ is an irreducible symplectic V-manifold, then $X$ has rational Gorenstein singularities and has a symplectic form on its smooth locus, hence by Theorem 6 of \cite{N3} it is a symplectic variety.

Moreover, if $f:\widetilde{X}\longrightarrow X$ is a resolution of the singularities, as $X$ has rational singularities we have an isomorphism between $H^{1}(X,\mathcal{O}_{X})$ and $H^{1}(\widetilde{X},\mathcal{O}_{\widetilde{X}})$. As $X$ has klt singularities we have $\pi_{1}(X)\simeq\pi_{1}(\widetilde{X})$ (by \cite{T}), and since $X$ is simply connected it follows that $\widetilde{X}$ is simply connected. As a consequence of this we see that $H^{1}(\widetilde{X},\mathcal{O}_{\widetilde{X}})=0$, and hence $H^{1}(X,\mathcal{O}_{X})=0$. As by definition of irreducible symplectic V-manifold we have that $H^{0}(X,\Omega_{X}^{[2]})$ is one dimensional, we are done.\endproof

No example of Namikawa symplectic variety which is not an irreducible symplectic variety, nor an irreducible symplectic V-manifold, nor a resolvable symplectic variety is known.

\subsubsection{Examples}

We now introduce two families of examples of irreducible symplectic varieties. In what follows $S$ will denote a projective K3 surface or an Abelian surface, and we let $\epsilon(S):=1$ if $S$ is K3, and $0$ if $S$ is Abelian. We let $\rho(S)$ be the rank of the N\'eron-Severi group $NS(S)$ of $S$.

An element $v\in \widetilde{H}(S,\mathbb{Z}):=H^{2*}(S,\mathbb{Z})$ will be written $v=(v_{0},v_{1},v_{2})$, where $v_{i}\in H^{2i}(S,\mathbb{Z})$, and $v_{0},v_{2}\in\mathbb{Z}$. It will be called \textit{Mukai vector} if $v_{0}\geq 0$, $v_{1}\in NS(S)$ and if $v_{0}=0$, then either $v_{1}$ is the first Chern class of an effective divisor, or $v_{1}=0$ and $v_{2}>0$. 

The $\mathbb{Z}-$module $\widetilde{H}(S,\mathbb{Z})$ has a pure weight-two Hodge structure and a lattice structure with respect to the Mukai pairing $(.,.)$ (see \cite{HL}, Definitions 6.1.5 and 6.1.11). We let $v^{2}:=(v,v)$ for every Mukai vector $v$, and we refer to $\widetilde{H}(S,\mathbb{Z})$ as the \textit{Mukai lattice} of $S$. We will always write $v=mw$, where $m\in\mathbb{N}$ and $w$ is a primitive Mukai vector on $S$. 

To any coherent sheaf $\mathcal{F}$ on $S$ we associate a Mukai vector $$v(\mathcal{F}):=ch(\mathcal{F})\sqrt{td(S)}\in\widetilde{H}(S,\mathbb{Z}).$$Taking $v$ a Mukai vector on $S$ and suppose that $H$ is a $v-$generic polarization (see as instance section 2.1 of \cite{PR3} for the precise definition), we consider the moduli space $M_{v}(S,H)$ (resp. $M_{v}^{s}(S,H)$) of Gieseker $H-$semistable (resp. $H-$stable) sheaves on $S$ with Mukai vector $v$. Then $M_{v}$ is a projective variety and $M_{v}^{s}\subseteq M_{v}$ is open (see \cite{HL}).

The following properties hold.
\begin{enumerate}
 \item If $S$ is a K3 surface and $v=mw$ for a primitive Mukai vector $w$, then $M_{v}\neq\emptyset$ if and only if $w^{2}\geq -2$ (see \cite{M2} and \cite{Y1}). If $S$ is an Abelian surface, then $M_{v}\neq\emptyset$ if and only $w^{2}\geq 0$ (see \cite{Y2}).
 \item If $S$ is a K3 surface and $v=mw$ for a primtive Mukai vector $w$ such that $w^{2}=-2$, then $M_{v}$ is a point (see \cite{M2}).
 \item If $S$ is a K3 surface, $v=mw$ and $w^{2}=0$, then $M_{v}\simeq Sym^{m}(S')$ for a projective K3 surface $S'$ (see \cite{M2} and \cite{KLS}). If $m=1$ we then get a projective K3 surface, while if $m\geq 2$ we then get an irreducible symplectic V-manifold, which is a resolvable symplectic variety but which is neither an irreducible symplectic variety nor an irreducible symplectic orbifold.
 \item If $S$ is an Abelian surface, $v=mw$ and $w^{2}=0$, then $M_{v}\simeq Sym^{m}(T)$ for an Abelian surface $T$ (see \cite{M2} and \cite{KLS}). If $m=1$ we then get an Abelian surface, and hence a K3 surface via the Kummer construction. If $m\geq 2$ we consider the natural sum morphism $s:M_{v}\longrightarrow T$, and let $K:=s^{-1}(0)$: then $K$ is an irreducible symplectic V-manifold, which is a resolvable symplectic variety but which is neither an irreducible symplectic variety nor an irreducible symplectic orbifold (see Example 1.13 of \cite{PR3}).
 \item If $S$ is a K3 surface or an Abelian surface and $v=mw$ for a primitive Mukai vector $w$ such that $w^{2}>0$, then $M_{v}$ is a normal, irreducible projective variety of dimension $v^{2}+2$ whose smooth locus is $M_{v}^{s}$ (see \cite{Y1}). By \cite{M1} $M_{v}$ has a symplectic form.
 \item If $S$ is an Abelian surface and $v=mw$ for a primitive Mukai vector $w$ such that $w^{2}>0$, by section 4.1 of \cite{Y2} we have a dominant isotrivial fibration $a_{v}:M_{v}(S,H)\longrightarrow S\times\widehat{S}$, where $\widehat{S}$ is the dual of $S$. We let $K_{v}:=a_{v}^{-1}(0_{S},\mathcal{O}_{S})$, and $K_{v}^{s}:=K_{v}\cap M^{s}_{v}$. The restriction of the symplectic form of $M_{v}$ to $K_{v}$ is a symplectic form (see \cite{Y2}).
\end{enumerate}

The moduli spaces $M_{v}$ and $K_{v}$ described above give us examples of irreducible symplectic varieties if $v^{2}>0$. This is the content of the following result, which is Theorem 1.19 of \cite{PR3}:

\begin{proposition}
\label{proposition:pr}
Let $S$ be a projective K3 surface of an Abelian surface, $v$ a Mukai vector on $S$ such that $v=mw$ for a primitive Mukai vector $w$ with $w^{2}>0$, and $H$ a $v-$generic polarization.
\begin{enumerate}
 \item If $S$ is K3, then $M_{v}(S,H)$ is an irreducible symplectic variety.
 \begin{enumerate}
  \item If $m=1$, it is an irreducible symplectic manifold.
	\item If $m=2$ and $w^{2}=2$, it has a symplectic resolution which is an irreducible symplectic manifold.
	\item In all other cases it has terminal singularities.
 \end{enumerate}
 \item Suppose that $S$ is Abelian. If $m=1$ and $w^{2}=2$ then $K_{v}(S,H)$ is a point. In all other cases $K_{v}(S,H)$ is an irreducible symplectic variety.
 \begin{enumerate}
  \item If $m=1$ and $w^{2}>2$, it is an irreducible symplectic manifold.
	\item If $m=2$ and $w^{2}=2$, it has a symplectic resolution which is an irreducible symplectic manifold.
	\item In all other cases it has terminal singularities.
 \end{enumerate}
\end{enumerate}
\end{proposition}

The case $m=1$ was proved in its final form by Yoshioka in \cite{Y1} and \cite{Y2} (but with important steps towards the complete proof given in \cite{M1}, \cite{B}, \cite{OG1}). The case $m=2$, $w^{2}=2$ was studied first by O'Grady in \cite{OG2} and \cite{OG3} for $v=2(1,0,-1)$, where it was shown that $M_{v}$ (resp. $K_{v}$) has a symplectic resolution which is an irreducible symplectic manifold. 

In \cite{LS} it is shown that the symplectic resolution exists for all $v=2w$, where $w$ is primitive and $w^{2}=2$. In \cite{PR1} it is shown that such a symplectic resolution is an irreducible symplectic manifold, deformation equivalent to $OG_{10}$ (resp. $OG_{6}$). The proof of the statement for all other cases is contained in \cite{PR3}. 

The cases $m=1$ and $m=2$, $w^{2}=2$ then recover all the known deformation classes of irreducible symplectic manifolds. Moreover, the case $m=2$, $w^{2}=2$ give examples of irreducible symplectic varieties which are resolvable symplectic varieties. We notice in particular that if $S$ is Abelian, then the smooth locus of $K_{v}$ in this case is not simply connected (see Theorem 3.7 of \cite{PR3}), giving then an example of irreducible symplectic variety which is not an irreducible symplectic orbifold.

The remaining cases give examples of irreducible symplectic varieties having no symplectic resolution, hence they are not resolvable symplectic varieties. Moreover, we have the following result:

\begin{proposition}
\label{proposition:defomk}
Let $S_{1}$ and $S_{2}$ be two projective K3 surfaces (resp. Abelian surfaces), $v_{i}=m_{i}w_{i}$ a Mukai vector on $S_{i}$ for $m_{i}>0$ and $w_{i}$ primitive Mukai vector, and $H_{i}$ a $v_{i}-$generic polarization on $S_{i}$. Then $M_{v_{1}}(S_{1},H_{1})$ (resp. $K_{v_{1}}(S_{1},H_{1})$) is deformation equivalent to $M_{v_{2}}(S_{2},H_{2})$ (resp. $K_{v_{2}}(S_{2},H_{2})$) if and only if $m_{1}=m_{2}$ and $w_{1}^{2}=w_{2}^{2}$.
\end{proposition}

The sufficient condition is essentially proved by Yoshioka in \cite{Y6} and \cite{Y7} (see even \cite{PR3}). The converse is proved in a forthcoming paper of the author and Rapagnetta (see \cite{PR4}). As a consequence we see that in dimension $2n$ we find a deformation class of irreducible symplectic varieties of dimension $2n$ for each pair $(m,k)$ such that $m^{2}k+1=n$ or $m^{2}k-1=n$ (just take the moduli space of semistable sheaves of Mukai vector $(m,0,-mk)$). 

By \cite{PR4} we see that the second Betti number of all these examples is 23 (in the case of $M_{v}(S,H)$ for a projective K3 surface $S$) or 7 (in the case of $K_{v}(S,H)$ for an Abelian surface $S$).

\begin{remark}
{\rm The deformation classes in dimension 4 and 6 we obtain in this way are different from those of the examples of singular irreducible symplectic V-manifolds we presented in section 2.1.2. Indeed, if $S$ is K3 we have $\dim(M_{v})=4,6$ only if $v^{2}=2,4$, so $v$ must be primitive and $M_{v}$ is then smooth. If $S$ is Abelian we have $\dim(K_{v})=4,6$ only if $v^{2}=6,8$. If $v^{2}=6$ then $v$ is primitive and $K_{v}$ smooth; if $v^{2}=8$ then either $v$ is primitive and $K_{v}$ is smooth, or $v=2w$ for $w$ primitive with $w^{2}$, and $K_{v}$ has a symplectic resolution.}
\end{remark}

Finally, we remark that the singular moduli spaces in the statement of Proposition \ref{proposition:pr} give examples of irreducible symplectic varieties which are not irreducible symplectic V-manifolds, as their singularities are not quotient singularities.


\section*{Acknowledgements}

The author wishes to thank Chiara Camere and Baohua Fu for useful discussions, and especially Gr\'egoire Menet and Antonio Rapagnetta for all the fruitful exchanges.

\end{document}